\documentclass[12pt]{amsart}
\usepackage{amssymb,amscd,diagrams,hyperref,ulem,color}

\newtheorem{thm}[subsection]{Theorem}

\newtheorem{lem}[subsection]{Lemma}

\newtheorem{prop}[subsection]{Proposition}
\newtheorem{cor}[subsection]{Corollary}

\newtheorem{rmk}[subsection]{Remark}
\newtheorem{ex}[subsection]{Example}
\theoremstyle{definition}
\newtheorem{defn}[subsection]{Definition}
\theoremstyle{remark}

\newenvironment{pf}{\paragraph{Proof}}{\par\medskip}
\numberwithin{equation}{subsection}

\newcommand\CC{\mathbb C}
\newcommand\ZZ{\mathbb Z}
\newcommand\QQ{\mathbb Q}

\newcommand{\op}{\mathrm}

\newcommand{\mc}{\mathcal}

\newcommand{\mb}{\mathbb}

\newcommand{\au}{\op{Aut}_0(S)}
\newcommand{\rest}[1]{{}_{\left|#1\right.}}
\begin{document}
\title[Automorphisms of    surfaces  with $\MakeLowercase{q}\geq 2$]{
Automorphisms of surfaces of general type with
$\MakeLowercase{q}\geq 2$  acting trivially in cohomology}

\author{Jin-Xing Cai}
\address{Jin-Xing Cai\\LMAM, School of  Mathematical Sciences\\Peking University\\Beijing 100871\\P.\,R.\,China}
\email{jxcai@math.pku.edu.cn}
\author{ Wenfei Liu}
\address{Wenfei Liu\\Fakult\"at f\"ur Mathematik\\Universt\"at Bielefeld\\Universit\"atsstr.\,25\\33615 Bielefeld\\Germany}
\email{liuwenfei@math.uni-bielefeld.de}
\author{Lei Zhang}
\address{Lei Zhang\\College of Mathematics and Information Sciences\\
Shaanxi Normal University\\Xi'an 710062\\P.\,R.\,China}
\email{lzhpkutju@gmail.com }

\thanks{The first author was supported by the  NSFC (No. 11071004). The second author was generously supported by the DFG via S\"onke Rollenske's Emmy-Noether Programm and by the Bielefelder Nachwuchsfonds.}

\begin{abstract}
In this note,  we prove that,  surfaces  of general type with
irregularity $q\geq 3$ are rationally cohomologically rigidified,
and so are minimal surfaces $S$ with $q(S)=2$  unless $K_S^2=8\chi(\mathcal O_S)$. Here a compact complex manifold
 $X$ is said to be rationally cohomologically rigidified  if its automorphism group
 $\op{Aut}(X)$ acts faithfully on
the cohomology ring $H^*(X, \mb Q)$.

As examples we give  a complete classification of surfaces isogenous
to a product with $q(S)=2$
  that are not rationally cohomologically rigidified.
\end{abstract}
\maketitle

\section{ Introduction}
 A compact complex  manifold $X$ is said to be  {\it
cohomologically rigidified} if its automorphism group $\op{Aut}(X)$
acts faithfully  on the cohomology ring $H^*(X,\ZZ)$, and  {\it
rationally cohomologically rigidified} if $\op{Aut}(X)$ acts
faithfully  on  $H^*(X,\QQ)$; it is
said to be  {\it  rigidified} if
 $\op{Aut}(X)\cap
\op{Diff}^0(X)=\{\op{id}_X\}$,  where $\op{Diff}^0(X)$ is the
connected component of the identity of the group of orientation
preserving diffeomorphisms of $X$ (\cite[Definition~12]{Cat11}).

 Note that any element in
$\op{Aut}(X)\cap \op{Diff}^0(X)$ acts trivially on the cohomology
ring $H^*(X,\ZZ)$. There are obvious implications: rationally
cohomologically rigidified $\Rightarrow$ cohomologically rigidified
$\Rightarrow$ rigidified.

  It is well known that curves of genus $\geq2$ are
   rationally cohomologically rigidified.  There are surfaces of general type
   with $p_g$ arbitrary large which are not cohomologically rigidified (\cite{Ca07}).
      An interesting
question posed by Catanese (\cite[Remark~46]{Cat11}) is whether every
surface of general type is rigidified.

The question  is closely related to the local moduli problem for
$X$, that is, whether  the natural local map
$\op{Def}(X)\to \mathcal{T}(M)_{[X]}$, from the Kuranishi space to the germ of the $\rm Teichm\ddot{u}ller$ space at $[X]$, is a homeomorphism or not. Here $M$ is
the underlying $C^\infty$-manifold of $X$ and  $[X]\in\mathcal{T}(M)$ is  the point
corresponding to the complex structure  of $X$ (\cite[Section
1.4]{Cat11}).

Apart from the local moduli problem, there is also motivation from
the global moduli problem, that is, the existence of fine moduli
space for polarized manifolds having the same Hilbert polynomial as
$X$ together with a so-called level $l$-structure (\cite[Lecture
10]{Po77}). Along this line, many authors  studied  the action of
automorphism groups of
    compact complex manifolds  on their cohomology rings.
  It is known that K3 surfaces are rationally cohomologically
rigidified (cf.\,\cite{BR75}, \cite{BHPV04}).
 For Enriques surfaces $S$,  either $S$ is cohomologically rigidified,
 or the kernel of $\op{Aut}(S)\rightarrow \op{Aut}(H^*(S,\ZZ))$
is a cyclic group of order $2$ or $4$, and the latter case was
completely  classified (\cite{Mu--Na84}, \cite{Mu10}). For elliptic
surfaces $S$, if $\chi (\mc O_S)>0$ and $p_g(S)>0$, then $S$ is
rationally cohomologically rigidified (\cite{Pet80}).
 There are also attempts at the generalization to hyperk\"ahler manifolds
        (\cite{Be83}, \cite{BNS11}); recently Oguiso \cite{O12}
        proved that  generalized Kummer manifolds are cohomologically rigidified.

 For surfaces of general
type, the problem seems harder,
 since there is not so many available structures on the cohomology groups
 as for the special surfaces. At the moment only partial results
are known.

Let $S$ be a  minimal nonsingular complex  projective  surface of
general type,  and $\op{Aut}_0(S)$   the subgroup of automorphisms
of $X$, acting trivially on  the cohomology ring $H^*(S, \mb Q)$.
Peters \cite{Pet79} showed that, if the canonical linear system
$\vert K_S\vert$ of $S$ is base-point-free, then $S$ is rationally
cohomologically rigidified, with the possible exceptional case where
$S$ satisfies either $K_S^2=8\chi(\mc O_S)$
 or $K_S^2=9\chi(\mc O_S)$.
In \cite{Ca06} and \cite{Ca10}, the first author proved that,
  if either $S$ has a  fibration of genus $2$ and
$\chi (\mathcal O_S)\geq 5$, or $S$ is an irregular surface   with
$K_S^2\leq 4\chi (\mathcal O_S)$ and  $\chi (\mathcal O_S)>12$, then
either $S$ is rationally cohomologically rigidified, or $\au$ is of
order two and   $S$ satisfies $ K_S^2=4\chi (\mathcal O_S)$ and $
q(S)=1$.

In this note,  we consider   surfaces of general type with $q(S)\geq
2$. Our main theorem is as follows.
\begin{thm}\label{main} Let $S$ be a  minimal nonsingular complex
   projective
surface of general type with $q(S)\geq 2$.  Then either $S$ is
rationally cohomologically rigidified, or $\au$ is a group  of order
two and $S$ satisfies $K_S^2=8\chi(\mathcal O_S)$, $q(S)=2$,    the
Albanese map of $S$ is surjective, and  $S$ has a pencil of genus
one.

In particular, if $q(S)\geq 3$, then $S$ is rationally
cohomologically rigidified, hence rigidified.
\end{thm}

By \cite[Theorem~45]{Cat11}, we have
\begin{cor}
Let $S$ be a minimal surface of general type with $q(S)\geq 3$. If
$K_S$ is ample, then the natrual map $Def(S)\rightarrow \mathcal
T(M)_{[S]}$ is a local homeomorphism between the Kuranishi space
and the Teichm\"uller space. Here $M$ is the underlying differential
manifold of $S$.
\end{cor}

 As examples we classify
surfaces isogenous to a product
  that are not rationally cohomologically rigidified.
\begin{thm}\label{example}
 Let $S=(C\times D)/G$ be a surface isogenous to a product with $q(S)\geq2$.
 Assume that  $S$ is not rationally cohomologically rigid.  Then $ S$ is as in Example~\ref{ex1} below; in particular, $S$ is of unmixed type,
 $G$ is isomorphic to one of the following groups:
 $\mathbb{Z}_{2m}\oplus\mathbb{Z}_{2mn}$,
 $\mathbb{Z}_{2}\oplus\mathbb{Z}_{2m}\oplus\mathbb{Z}_{2mn}$ ($m, n$ are
 arbitrary  positive integers), and
$g(C/G)=g(D/G)=1$.
\end{thm}

By a result of Borel-Narasimhan (\cite{BN67}), surfaces in Example~\ref{ex1} are rigidified (Proposition~\ref{bn}).  It is not known
whether  surfaces $S$ in the latter case of  Theorem~\ref{main} are
rigidified. A further step that can be done on these surfaces is to
check if the action of $\op{Aut}(S)$ on cohomology with
$\ZZ$-coefficient or on the fundamental group is faithful.

The paper is organized as follows.

Theorems~\ref{main} and
\ref{example} are proved in Sections 2--3 and Section 4,
respectively.

In Section 2, we consider projective manifolds (of arbitrary
dimension) with maximal Albanese dimension. We show that the group
of automorphisms of such varieties behave like that of curves
(Theorem~\ref{mad}), which is of independent interest.

In Section~3, a combination of the topological and holomorphic
Lefschetz formulae for a group action and the Severi inequality
helps us pin down the numerical restrictions on the surfaces with
non-trivial $\op{Aut}_0(S)$.

In the classification of surfaces isogenous to a product $S=(C\times
D)/G$ that are not rationally cohomologically rigidified in
Section~4, we first use the theorem proved in Section~2 to exclude
the mixed type case; then Broughton's cohomology representation
theorem for curves is used to calculate the cohomology of surfaces
isogenous to a product of unmixed type; finally we manage to give
the classification by finding an appropriate character of the group
$G$ through Frobenius' reciprocity theorem.

\

\noindent{\bf Notations.} A pencil of genus $b$  of a surface $S$ is
a fibration $f:S\to B$, where  $B$ is a smooth curve of genus $b$.

For a  smooth projective variety $X$,  we denote by $p_g(X)$,
$q(X)$, $e(X)$, $\chi(\mathcal{O}_X)$, $\chi(\omega_X)$, and $K_X$
the geometric genus, the irregularity, the Euler topological
characteristic,  the Euler characteristic of the structure sheaf, the
Euler characteristic of the canonical  sheaf, and a canonical divisor of
$X$, respectively.

We denote by $\op{Aut}_0(X)$ the kernel the natural homomorphism of
groups $\op{Aut}(X)\rightarrow \op{Aut}(H^*(X,\QQ))$.

 We use $\equiv$,
$\equiv_\mathbb{Q}$ to denote linear equivalence and
$\mathbb{Q}$-linear equivalence of divisors, respectively.

For a finite group $G$ and an element $g\in G$,  we denote by

 $|G|$ : the order of $G$,

$|g|$ : the order of  $g$,

 $C_G(g)$ : the conjugacy class of $g$ in $G$,

 $\op{Irr}(G)$ : the
set of irreducible characters of $G$,

$\op{Ker}(\chi):=\{g\in G \mid \chi(g)=\chi(1)\}$, for $\chi\in
\op{Irr}(G) $.

 For a representation $V$
of $G$ and a character $\chi\in \op{Irr}(G)$, we let $V^\chi$ be the
sum of irreducible sub-$G$-modules $W$ of $V$ with $\chi_W=\chi$,
where  $\chi_W$ is the character of $G$-module $W$.

Let  $H$ be a subgroup of a finite group $G$,
 and $\chi$  a character of $H$,
 we denote by $\chi^G$  the induced character from $\chi$.
Recall that  $\chi^G$ is defined by
\[
 \chi^G(g) = \frac{1}{|H|}\sum_{t\in G} \chi^\circ (tgt^{-1})
\]
where for any $g\in G$
\[
 \chi^\circ (g) =
\begin{cases}
 \chi(g), & \text{ if } g\in\, H, \\
 0, & \text{ if } g\notin\, H.
\end{cases}
\]

The symbol $\mathbb{Z}_n$ denotes  the cyclic group of order $n$.

\

{\bf Acknowledgements.} Part of this note appears in the
doctoral thesis of the second author submitted to Universit\"at
Bayreuth (2010). The second author and the third author would like to
thank the hospitality of Professor Jin-Xing Cai when they visited
School of Mathematical Sciences, Beijing University in August 2012.

\section{Projective manifolds with maximal Albanese dimension}
In this section, we use the generic vanishing theorems of Green and
Lazarsfeld (\cite{GL87}, \cite{GL91}, see also \cite{Ha04}) and the
notion of continuous global generation  (\cite{PP03}, \cite{PP06},
\cite{BLNP12}) to show that the groups of automorphisms of projective
manifolds of general type with maximal Albanese dimension  and with
positive generic vanishing index behave like that of curves
  (Theorem~\ref{mad}).

We begin by recalling some notations.
\subsection{The generic vanishing index} Let $X$ be a  smooth projective variety with $q(X)>0$, and
$a: X\to \op{Alb} X$ the Albanese map of $X$. We say that $X$ is of
maximal Albanese dimension if $a$ is a generically finite map onto
its image.

 For $0 \leq i\leq \dim X$, the  $i$-th cohomological
support locus of  $X$ is defined as
$$V^i(\omega_X) := \{\alpha \in \op{Pic}^0(X)\,|\,
 h^i(X, \omega_X \otimes \alpha) > 0\}.$$
Let $$\op{gv}_i(\omega_X)=
  \op{codim}_{\op{Pic}^0(X)}V^i(\omega_X)-i,\text{ and }\
\op{gv}(\omega_X) = min_{i>0}\{ \op{gv}_i(\omega_X)\}.$$ We call
$\op{gv}_i(\omega_X)$ and $\op{gv}(\omega_X)$  the $i$-th generic
vanishing index and the generic vanishing index of $X$,
respectively.
\subsection{The results of Green and Lazarsfeld}
\label{gl}If $X$ is a   smooth projective variety of
maximal Albanese dimension, by  the generic vanishing theorem due to
Green and Lazarsfeld (cf. \cite{GL87}, \cite{GL91}),  one has
\begin{itemize}\item[(\ref{gl}.1)]
$\op{gv}(\omega_X)\geq0$. So for a general $\alpha\in
\op{Pic}^0(X)$, $h^i(X, \omega_X\otimes\alpha)=0$ for all $i>0$, and
hence  $\chi( \omega_X)=\chi( \omega_X\otimes\alpha)=h^0(X,
\omega_X\otimes\alpha)$;
\item[(\ref{gl}.2)]
for each positive dimensional  component $Z$ of $V^i(\omega_X)$, $Z$
is a complex sub-torus of $\op{Pic}^0(X)$, and there exists an
algebraic variety $Y$ of dimension $\leq \dim X-i$ and a dominant
map $f: X \rightarrow Y$ such that $Z \subset \alpha +
f^*\op{Pic}^0(Y)$ for some  $\alpha\in \op{Pic}^0(X)$.
\end{itemize}

\subsection{A result on varieties with positive generic vanishing indices and its applications}
\begin{thm}\label{euln}
Let $f:X\to Y$ be a generically finite morphism of  smooth
projective varieties of maximal Albanese dimension. Let $a: X\to
A:=\op{Alb} X$ be the Albanese map of $X$. Assume  that   the
following conditions hold:
\begin{itemize}
\item[(i)]$a$ factors through $f$;
\item[(ii)]$\op{gv}_i(\omega_X)\geq 1$ for all $0<i<\dim X$;
\item[(iii)] $p_g(X) = p_g(Y)$ if $q(X)=\dim X$.
\end{itemize}
Then  $\chi(\omega_X) \geq \chi(\omega_Y)$ and $``="$ occurs only
when $f$ is birational.
\end{thm}
\begin{pf}
By the assumption (i), there is a morphism $a': Y\to A$, such that
$a'\circ f=a$. By the universal property of the Albanese map, we
have that $a'$ is just the Albanese map of $Y$.

 Since $f$ is generically finite and $Y$ is smooth, there is  an
injective   morphism of sheaves $ f^*\omega_Y \hookrightarrow
\omega_{X}$. Taking $f_*$ and composing with the natural morphism
$\omega_Y \rightarrow f_*f^*\omega_Y$, we obtain an inclusion of
sheaves $\omega_Y \hookrightarrow  f_*\omega_{X}$.

Taking $a'_*$, we get an inclusion of sheaves
$$\rho: a'_*\omega_Y \hookrightarrow  a'_*(f_*\omega_X)=a_*\omega_X.$$
Hence for every $\alpha \in \op{Pic}^0 A$, we have   an inclusion
$$\rho_\alpha: H^0(A, a'_*\omega_Y \otimes \alpha)
 \hookrightarrow  H^0(A, a_*\omega_X\otimes \alpha).$$
By (\ref{gl}.1), we have $\chi(\omega_X) \geq \chi(\omega_Y)$ by
choosing $\alpha$ to be general.

Now  we will show that, if $\chi(\omega_X) = \chi(\omega_Y)$ then
$\deg f=1$.

Again by (\ref{gl}.1), the assumption $\chi(\omega_X) =
\chi(\omega_Y)$
 implies that,  for a general $\alpha\in
\op{Pic}^0(X)$, $h^0(X, \omega_X\otimes\alpha)=h^0(Y,
\omega_Y\otimes\alpha)$. Thus  we can find a non-empty Zariski open
set $U \subset \op{Pic}^0(A)$ such that for $\alpha \in U$,
$\rho_\alpha$ is an isomorphism.  Consider  the following
commutative diagram
\begin{equation}\label{eq: ccg}
 \begin{CD}
\bigoplus_{\alpha \in T}H^0(A, a'_*\omega_Y \otimes \alpha)
\otimes \alpha^{-1} @> \op{ev}_T'>>a'_*\omega_Y\\
@V\oplus_{\alpha \in T}\rho_\alpha \otimes \alpha^{-1} VV
 @V\rho VV   \\
\bigoplus_{\alpha \in T}H^0(A, a_*\omega_X\otimes \alpha)\otimes
\alpha^{-1}   @>\op{ev}_T >>   a_*\omega_X
\end{CD}
\end{equation}
where $T$ is a subset of $\op{Pic}^0(A)$ and
  $\op{ev}_T$, $\op{ev}_T'$ are
 evaluation maps.

{\it Case 1. $q(X)>\dim X$}.  We let $T=U$. Then by the choice of
$U$, $\oplus_{\alpha \in T}\rho_\alpha$ is an isomorphism.
 In this case  the
 assumption ($ii$)  is equivalent to
$\op{gv}(\omega_X)\geq 1$.  This implies that $\op{ev}_T$ is
surjective  (\cite[Proposition~5.5]{PP06}, \cite[(a) of Corollary
4.11]{BLNP12}).  By the commutative diagram \eqref{eq: ccg},  it
follows that $\rho$ is surjective.

{\it Case 2. $q(X)=\dim X$}. In this case we let $T=U\cup
\{\hat{0}\}$, where $\hat{0}$ is the identity element of
$\op{Pic}^0(A)$.  By the assumption ($iii$)  and by the choice of
$U$, we have that $\oplus_{\alpha \in T}\rho_\alpha\otimes
\alpha^{-1}$ is an isomorphism. By \cite[(b) of Corollary
4.11]{BLNP12},
 the  assumption ($ii$) implies  $\op{ev}_T$ is surjective, and
  so $\rho$ is surjective.

Since the ranks (at the generic point of $A$) of $a'_*\omega_Y$,
$a_*\omega_X$ are $\deg a'$, $\deg a$ ($=\deg a'\cdot\deg f$),
respectively,   the surjection of $\rho$  implies $\deg f=1$.

 This  completes the  proof of the theorem.
\qed\end{pf}

\begin{rmk}
After finishing the paper, we are kindly informed by Sofia Tirabassi that, under slightly milder hypothesis, Theorem 2.4 was already independently  proved in her thesis (\cite[Proposition 5.2.4]{T11}) in the case where $q(X) > \dim X$.
\end{rmk}

\begin{thm}\label{mad}
Let $X$ be a   smooth projective variety of general type and  of
maximal Albanese dimension.  If
      $\op{gv}_i(\omega_X)\geq 1$ for all $0<i<\dim X$,
then   $X$ is  rationally cohomologically rigidified.
\end{thm}
\begin{pf}Otherwise, there is a non-trivial automorphism $\sigma$  of $X$,
which acts trivially on $H^*(X, \mathbb Q)$. Since  $X$ is of
general type,
   $\sigma$  is of  finite order. Replacing $\sigma$ by a suitable power, we
   may assume $\sigma$ is of  prime order, say $p$.

Let $$\pi:X \rightarrow \bar X= X/\left<\sigma\right>$$ be the
quotient map. Since  $\sigma$ acts trivially on $ H^i(X,
\mathbb{C})$ for all $i\geq 0$, by Hodge theory, we have
$$H^i(\bar X, \mathcal{O}_{\bar X})\simeq H^i(X, \pi^\sigma_*\mathcal{O}_X)
 =H^i(X,
\mathcal{O}_X)^\sigma = H^i(X, \mathcal{O}_X).$$ In particular, we
have
\begin{align}\label{2.4.1}
h^{\dim X}(\bar X, \mathcal{O}_{\bar X}) = h^{\dim X}(X, \mathcal{O}_X), \textrm{
and }  \chi(\mathcal{O}_{\bar X}) = \chi(\mathcal{O}_X).
\end{align}
 Let $\rho: Y \rightarrow \bar X$ be a resolution of quotient
singularities (if any). Then $R^i\rho_*\mathcal{O}_Y = 0$ for $i>0$
since  quotient singularities  are rational. Thus
\begin{align}\label{2.4.2}h^{\dim X}(Y,
\mathcal{O}_Y) = h^{\dim X}(\bar X, \mathcal{O}_{\bar X}), \textrm{ and }
\chi(\mathcal {O}_{Y}) = \chi(\mathcal{O}_{\bar X}).
\end{align}

By (\ref{2.4.1}) and  (\ref{2.4.2}), using Serre duality, we obtain
$$p_g(Y) = p_g(X), \textrm{
and }  \chi(\omega_{Y}) = \chi(\omega_X).$$

We claim that $X^\sigma\not=\emptyset$. Otherwise, the  map $\pi$ is
$\rm\acute{e}$tale. This implies $\chi(\omega_X) =
p\chi(\omega_{\bar{X}})$. Combining this with (\ref{2.4.1}), we have
$\chi(\omega_X)=0$.  On the other hand, since
$\op{gv}_i(\omega_X)\geq 1$ for all $0<i<\dim X$ by assumption,
$\chi(\omega_X)=0$ is equivalent to   $X$ being not of general type
(\cite[Proposition~4.10]{BLNP12}). So we get a contradiction.

Let $a: X\to \op{Alb} X$ be the Albanese map of $X$ (the map $a$ is
unique up to translations of $\op{Alb} X$  and we fix it once for
all). We have that there is an automorphism $\bar \sigma$ of
$\op{Alb} X$, such that $\bar \sigma\circ a=a\circ \sigma$. Since
$\sigma$ induces
 trivial action  on  $H^1(X, \mathbb Q)$, we have  that
 either $\bar \sigma$ is a translation or $\bar \sigma=\op{id}_{\op{Alb} X}$.
 If $\bar \sigma$ is a translation, then $X^\sigma=\emptyset$ --- a contradiction by the claim above.  So $\bar
\sigma=\op{id}_{\op{Alb} X}$, and consequently,
   $a$ factors
through $\pi$.

 Let   $f:X\to Y$ be the rational  map induced by the quotient map $\pi$, and
   $\rho: X'\to X$ be a birational morphism such that $f\circ\rho$
is a morphism.
 Since  $\op{V}^i(\omega_X)$ are birational invariants,
 using $X'$ instead of $X$ and $a\circ\rho$ instead of $a$,  we may
 assume that $f$ is a morphism.
  Then $\deg f$ is divisible by $p$, which is a contradiction
    by Theorem~\ref{euln}.
 \qed\end{pf}
\begin{cor}\label{c1}
Let $S$ be a   smooth projective surface of maximal Albanese
dimension.
    Assume  that  $S$ has
no pencils of  genus  $ \geq 2$, and  $S$ has no pencils of  genus $
1$ when $q(S) = 2$. Then $S$ is  rationally cohomologically
rigidified.
\end{cor}
\begin{pf}
  Note that  the assumption  of Corollary
\ref{c1} is equivalent to $\op{gv}_1(\omega_S)\geq1$ by
(\ref{gl}.2). The corollary follows by Theorem~\ref{mad}.
 \qed\end{pf}
\begin{rmk}\rm   The assumption  on the  $i$-th generic vanishing index
(which is slightly weaker than $\op{gv}(\omega_X)\geq1$)
 is indispensable for Theorem~\ref{mad}.
For example, let  $X=S\times C$, $\tau=\sigma\times\op{id}_C$, where
the pair ($S, \sigma$) is as in Example~\ref{ex1}, and $C$ is a
curve of genus $\geq2$. Then $\tau$ is an involution of $X$, which
acts trivially on  the cohomology ring $H^*(X, \mathbb Q)$.

 It is interesting to  classify   smooth projective $3$-folds of
general type and of  maximal Albanese dimension that are not
rationally cohomologically rigidified.
\end{rmk}

\section{Numerical classifications}
\begin{thm}\label{numer} Let $S$ be a  minimal nonsingular complex
   projective
surface of general type with $q(S)\geq 2$. Then either  $\au$ is
trivial, or $\au$ is a group  of order two and $S$ satisfies
$K_S^2=8\chi(\mathcal O_S)$.
\end{thm}
 For the proof of Theorem~\ref{numer}, we  need the following lemmas.
\begin{lem}\label{x1}
Let $S$ be a  minimal nonsingular complex    projective surface of
general type, and $G\subset \op{Aut}(S)$ be a subgroup of
automorphisms of $S$. Assume that the quotient $S/G$ is of
 general type. Then $K_S^2\geq |G|K_X^2$, where
  $X$ is the    minimal smooth model of  $S/G$.
\end{lem}
\begin{pf}
Let $Y=S/G$, and $\pi: S\to Y$ be the quotient map.
 Let $\eta:
\tilde{Y}\to Y$ be the minimal resolution of the quotient
singularities of $Y$.  We have  a commutative diagram
$$\begin{diagram}
  \tilde S  &\rTo^{\tilde \pi} & \tilde{Y} \\
  \dTo^\alpha   &  &  \dTo^\eta    \\
S & \rTo^\pi  & Y
\end{diagram}$$
 where  $\tilde S$ is  the minimal resolution of the normalization of
  the fiber product $S\times_Y \tilde{Y}$, and $\alpha$,
  $\tilde{\pi}$ are natural induced maps.

Since  surface quotient singularities  are log-terminal (cf.\,\cite[Corollary~1.9]{Ka84}), we have that  there is an effective
$\mathbb{Q}$-divisor $D$ supported on $\eta$-exceptional curves,
such that
\begin{equation}\label{1.1}K_{\tilde Y}\equiv_\mathbb{Q}\eta^*K_{Y}-
D.\end{equation}
 Let $R$  be the ramification divisor  of $\pi$, that is, $
 R=\sum_C{(|G_C|-1)C}$
   with the sum taken all smooth curves $C$ on $S$
 and   $G_C=\{g\in G\ |\ g\rest
{C}=\op{id}\}$. We have
\begin{equation}\label{1.2}K_{ S}\equiv_\mathbb{Q}\pi^*K_{Y}+R.
\end{equation}
Let $\rho:\tilde{Y}\to X$ be  the contraction of $\tilde Y$ onto its
uniquely determined minimal model $X$. Then
\begin{equation}\label{1.3}K_{\tilde Y}=\rho^*K_{X}+A\end{equation}
 for some effective exceptional divisor $A$ of
$\rho$. Note
that $\pi\circ\alpha=\eta\circ\tilde{\pi}$, and we have
\begin{align}\label{can11}\alpha^*K_{S}\equiv_\mathbb{Q}
(\rho\circ\tilde{\pi})^*K_{X}+
\tilde{\pi}^*A+\tilde{\pi}^*D+\alpha^*R
\end{align}
by combining (\ref{1.1}), (\ref{1.2}) and (\ref{1.3}).
So
\begin{align*} K_{S}^2=&\alpha^*K_{S}((\rho\circ\tilde{\pi})^*K_{X}+
\tilde{\pi}^*A+\tilde{\pi}^*D+\alpha^*R)\\
\geq&\alpha^*K_{S}(\rho\circ\tilde{\pi})^*K_{X}
\ \ \ ( \textrm{ since $\alpha^*K_{S}$ is nef}) \\
\geq&((\rho\circ\tilde{\pi})^*K_{X})^2 \ \ \
 (\textrm{ using (\ref{can11}), since $(\rho\circ\tilde{\pi})^*K_{X}$ is nef)} \\
 =&|G|K_X^2.
\end{align*}
\qed\end{pf}
\begin{lem} \label{3.1}
Let $S$ be a  minimal nonsingular complex
   projective
surface of general type with $q(S)\geq 2$. If $S$ has a pencil of
genus larger than one, then $\au$ is trivial.\end{lem}
\begin{pf}
 Let  $f\colon S\to
B$ be such a fibration  over a curve $B$ of genus $b\geq2$.  Suppose
that there is  a non-trivial element  $\sigma\in \au$.
   Let $F$ be a
 general fiber of $f$, and  $g$  the genus of $F$. We have $g\geq2$.
  Since $\sigma$
acts trivially on $\op{NS}(S)\otimes\Bbb Q\hookrightarrow H^2(S,
\Bbb Q)$, we have  $\sigma^*F$ is numerically equivalent to $F$. So
$f$ is preserved under the action of $\sigma$. Since $\sigma$ acts
trivially on $f^*H^1(B, \mathbb{Q})\subseteq H^1(S, \mathbb{Q})$,
 it induces identity  action on $B$, and so
     $f\circ\sigma=f$.  We have $b\leq1$ by \cite[Lemma~2.1]{q2} --- a
     contradiction.
\qed\end{pf}

\subsection{Proof of Theorem~\ref{numer}}
 By Lemma~\ref{3.1},  we may assume that $S$ is of
maximal Albanese dimension.

Assume that $G:=\au$ is not trivial.  We will show that $G$ is a
group of order two and $S$ satisfies $K_S^2=8\chi(\mathcal O_S)$.

  Let $X$ be a minimal  smooth model of the
quotient $S/G$. Then $q(X)=q(S)$, $p_g(X)=p_g(S)$, and  both the
canonical map and the Albanese map of $S$ factorize  through the
quotient map $ S\to S/G$. In particular, $X$ is of maximal Albanese
dimension and
 the Kodaira
 dimension of $X$ is at least $1$.

   If the Kodaira
 dimension of $X$ is  $1$, then  the canonical map $\phi_X$ is
 composed with a pencil of genus $= q(X)-1$ (cf.\,\cite[p.\,345, Lemme]{Be82}).
Since $\phi_S$ factors through the quotient map $ S\to S/G$,
    we have that
   $\phi_S$ is composed with a pencil  whose base curve, say $C$,
    is
   of genus $g(C)\geq q(X)-1=q(S)-1\geq1$.  On the other hand,
    one has that
either $q(S)=g(C)=1$ or $g(C)=0$ and $q(S)\leq2$ by
\cite{Xi85}.
    This is a contradiction.

 Now we may assume  that $X$
is of
 general type.
 By Severi inequality (\cite{Par05}), we have
\begin{equation}\label{se}K_X^2\geq4\chi(\mathcal
O_X)=4\chi(\mathcal O_S).\end{equation}
 Combining
(\ref{se}), Lemma~\ref{x1} and the Bogomolov--Miyaoka--Yau inequality
$9\chi(\mathcal O_S)\geq K_S^2$, we get $|G|=2$ and
\begin{equation}\label{e1} K_S^2\geq8\chi(\mathcal O_S).
\end{equation}

Let $\sigma$ be the generator of $G$. Let $D_i$ ( $1\leq i\leq u$,
$u\geq 0$) be the $\sigma $-fixed curves. After suitable
re-indexing, we may assume that $D_i^2\geq0$ for $i\leq k$ ($0\leq
k\leq u$) and $D_i^2<0$ for $i>k$. We may apply
 the topological  and holomorphic Lefschetz formula to $\sigma $ (cf.\,[AS, p.\,566]):
\begin{equation}\label{eq: lefschetz}
\begin{split}
&e (S)+8(q(S)-\op{dim}_\CC H^0(S, \Omega _S^1)^\sigma )-
2(h^2(S, \mb C)- \op{dim}_\CC H^2(S, \mb C)^\sigma  )\\
&=e (S^\sigma )=n+\sum_{i=1}^u e(D_i)\\
&\chi (\mc O_S)+2(q(S)-\op{dim}_\CC H^0(S, \Omega _S^1)^\sigma )
={n-\sum_{i=1}^u K_SD_i\over 4},
\end{split}
\end{equation}
where $n$ is the number of isolated $\sigma $-fixed points.
Combining \eqref{eq: lefschetz} with Noether's formula, we get
\begin{equation}\label{top}
K_S^2=8\chi(\mc O_S)+\sum_{i=1}^u D_i^2\leq8\chi(\mc
O_S)+\sum_{i=1}^k D_i^2.
\end{equation}

Let  $\rho:\tilde S\to S $ be the blowup of all
 isolated fixed points of $\sigma $, and $\tilde \sigma $ the induced involution
 on $\tilde S$. Let $\tilde \pi: \tilde S\to\tilde{X}:=\tilde S/\tilde \sigma$
be the quotient map. Let $\eta:\tilde{X}\to X$ be  the map
contracting all $(-1)$-curves on $\tilde{X}$. We have
\begin{equation}\label{can}\rho^*K_{S}=(\eta\circ\tilde{\pi})^*K_{X}+
\tilde{\pi}^*A+\sum_{i=1}^u\rho^*{D_i}
\end{equation} for some effective exceptional
divisor $A$ of $\eta$.

We show  that  $k=0$. Otherwise,  we have
\begin{align*} K_{S}^2=&\rho^*K_{S}^2\\
\geq&\rho^*K_{S}(\eta\circ\tilde{\pi})^*K_{X}+
\sum_{i=1}^k\rho^*K_{S}\rho^*{D_i}
\ \ \ (\textrm{using (\ref{can}), since $\rho^*K_{S}$ is nef}) \\
\geq&((\eta\circ\tilde{\pi})^*K_{X})^2+\sum_{i=1}^k\rho^*K_{S}\rho^*{D_i}
\ \ \
 (\textrm{ using (\ref{can}), since $(\eta\circ\tilde{\pi})^*K_{X}$ is nef)} \\
 =&2K_X^2+\sum_{i=1}^kK_{S}D_i\\
\geq&8\chi(\mathcal
O_S)+(K_{S}-\sum_{i=1}^kD_i)\sum_{i=1}^kD_i+\sum_{i=1}^kD_i^2 \ \ \
(\textrm{by
(\ref{se})})\\
 \geq& 8\chi(\mathcal O_S)+2+\sum_{i=1}^kD_i^2,\\
\end{align*}
which contradicts (\ref{top}), where the last inequality follows
since  each $\sigma$-fixed curve is contained in the fixed part
 of $|K_S|$ (cf.\,\cite[1.14]{Ca04}) and $|K_S|$ is $2$-connected
 (cf.\,\cite[VII, Proposition~6.2]{BHPV04}).

 So we have $k=0$ and hence $u=0$ by combining
(\ref{e1}) with (\ref{top}). This finishes the proof of Theorem~\ref{numer}.\qed
\subsection{Proof of Theorem~\ref{main}}
By Theorem~\ref{numer}, there remains to prove the following claim:
if $\au$ is not trivial, then $q(S)=2$,  the Albanese map of $S$ is
surjective, and $S$ has a pencil of genus one.

By Lemma~\ref{3.1}, we may assume that $S$ is of maximal Albanese
dimension, and $S$ has no  pencils of genus $\geq2$. Now the claim
follows by Corollary \ref{c1}. \qed
\section{Surfaces isogenous to a product}
Surfaces isogenous to a product  play an important role in studying
the geometry and the moduli of surfaces of general type. For
examples, they provide simple counterexamples to the DEF=DIFF
question whether deformation type and diffeomorphism type coincide
for algebraic surfaces (\cite{Cat03}), and they are useful in the
construction of Inoue-type manifolds (\cite[Definition~0.2]{BC12}).
In this section we give an explicit description for surfaces $S$
isogenous to a product  with $q(S)\geq2$ which are not rationally
cohomologically rigidified (Examples~\ref{ex1} and Theorem~\ref{non-abelian}).

We begin by recalling some notations of surfaces isogenous to a
product; we refer to \cite{Cat00}
  for properties of these surfaces.

\begin{defn}(\cite[Definition~3.1]{Cat00}) A smooth projective surface
 $S$ is isogenous to a (higher) product if it is a quotient $S=(C\times D)/G$,
 where $C,D$ are curves of genus at least two,
 and $G$ is a finite group acting freely on $C\times D$.
\end{defn}

 Let $S=(C\times D)/G$ be a surface isogenous to a product.
 Let
$G^\circ$ be the intersection of $G$ and $\op{Aut} (C)\times
\op{Aut}(D)$ in $\op{ Aut}(C\times D)$. Then $G^\circ$ acts on the
two factors $C,D$ and  acts on $C\times D$ via the diagonal action.
If $G^\circ$ acts faithfully on both $C$ and $D$, we say $(C\times
D)/G$ is a minimal realization of $S$. By \cite[Proposition~3.13]{Cat00}, a minimal realization exists and is unique. In the
following we always assume $S=(C\times D)/G$ is the minimal
realization.

We say that $S$ is of unmixed type if $G=G^\circ$;
  otherwise $S$ is of mixed type.
\begin{prop}\label{mtype}
 If $S=(C\times C)/G$ is a surface isogenous to a product  of mixed type with
 $q(S)\geq2$,
  then
  $S$ is  rationally cohomologically rigidified.
\end{prop}
 \begin{pf}
Let $\sigma\in G\setminus G^\circ$. Up to coordinate change in both factors
of $C\times C$, we can assume
$\sigma(x,y) = (y,\tau x)$ for some $\tau\in G^\circ$
(cf.\,\cite[Proposition~3.16]{Cat00}).
 Then $X := (C \times C)/\sigma$ is smooth, and the natural maps
  $C\times C\rightarrow X$ and $\pi\colon X\rightarrow S$ are both \'etale coverings.

Note that
\begin{equation}
\begin{split}
&\op{Pic}^0(X) \cong \op{Pic}^0(C \times C)^{\sigma}\cong (\op{Pic}^0(C) \times \op{Pic}^0(C ))^{\sigma}  \\
&=\{(\alpha, \beta) \in  \op{Pic}^0(C) \times \op{Pic}^0(C )|
\alpha= \tau^*\beta,\beta = \alpha\}
\end{split}
\end{equation}
Hence we can identify $\op{Pic}^0(X)$ with the set $$\{(\alpha,
\alpha)\in  \op{Pic}^0(C) \times \op{Pic}^0(C )| \tau^*\alpha =
\alpha\}.$$ We have
\begin{equation}
\begin{split}
&H^1(X, (\alpha, \alpha)) \\
&\cong H^1(C\times C, (\alpha, \alpha))^{\sigma}\\
&\cong (H^1(C, \alpha) \otimes_\CC H^0(C, \alpha) \oplus H^0(C, \alpha) \otimes_\CC H^1(C, \alpha))^{\sigma}\\
\end{split}
\end{equation}
which is zero unless $\alpha = \hat{0}$, the identity element of
$\op{Pic}^0(C)$. Using Serre duality we have $V^1(\omega_X) =
\{\hat{0}\}$.

Note that  $\pi\colon X\rightarrow S$ is an \'etale covering. We
have for any $\gamma\in \op{Pic}^0(S)$
\begin{equation}\label{eq: gvs}
  H^i(X, \omega_X \otimes \pi^*\gamma) = H^i(S, \omega_S\otimes \gamma \otimes \pi_* \mc O_X)
\end{equation}
by the projection formula and the Leray spectral sequecne. The left hand
side of \eqref{eq: gvs} is zero unless $\pi^*\gamma = \hat 0$, while
the right hand side contain a direct summand of $H^i(S,
\omega_S\otimes \gamma)$. In other words, $H^i(S, \omega_S\otimes
\gamma) =0 $ unless $\pi^*\gamma = \hat 0$. Since  $\pi^*\colon
\op{Pic}^0(S) \rightarrow \op{Pic}^0(X)$ is a finite map onto its
image, we conclude that $V^1(\omega_S)$ is a finite set. In
particular, we have $\op{gv}_1(\omega_S)~\geq~1$. By Theorem~\ref{mad}, we have that $S$ is  rationally cohomologically
rigidified. \qed\end{pf}

Contrary to the case of surfaces isogenous to a product of mixed
type, there are surfaces isogenous to a product of unmixed
type which are  not rationally cohomologically rigidified. Before
giving such examples,  we insert here two facts on curves as well as an expression for the second cohomology of surfaces isogenous to a product of unmixed type that will
be used in the sequel.

\subsection{Riemann's existence theorem}\label{riemann}
 Let  $m_1,\dots,m_r\geq2$ be $r$
 integers, and  $G$  a finite group. Let
 $B$ be a curve of genus $b$, and let $p_1, \cdots, p_{r}\in B$
  be $r$ different points.

Assume that there are $2b+r$ elements of $G$ (not necessarily
different),  $\alpha_j, \beta_j, \gamma_i$ ($1\leq j\leq b$, $1\leq
i\leq r$), such that   these elements generate $G$, and satisfy
\begin{equation}\label{rie1}
\prod_{j=1}^b[\alpha_j\beta_j]\prod_{i=1}^r\gamma_i=1,\ \textrm{and
} \ |\gamma_i|=m_i.
\end{equation}

If the Riemann--Hurwitz equation
\begin{equation*}\label{rie} 2g-2=|G|(2b-2+\sum_{i=1}^r(1-{1\over
m_i}))
\end{equation*}
 is satisfied for some non-negative integer $g$,
then
  there exists a curve $C$ of
genus $g$
  with a faithful $G$-action, such that the
quotient map $C\to C/G\simeq B$ branched  exactly at $p_1, \cdots,
p_{r}$, and $m_i$ is
  the ramification numbers of  $q$ over
$p_i$.

In what follows  we call a $2b+r$-tuple ($\alpha_1,\cdots, \alpha_b,
\beta_1, \cdots,  \beta_b, \gamma_1, \cdots, \gamma_i$) of elements
of $G$ a generating vector  of type $(b; m_1,\dots,m_r)$, if these
$2b+r$ elements generate $G$ and satisfy (\ref{rie1}).

\subsection{The cohomology representation of the   group of automorphisms of
a curve}
 Let $C$ be a smooth curve of genus $g(C) \geq 2$
 and $G$ a   group of automorphisms of $C$.
Let $r$ be the number of branch points of the quotient map $
C\rightarrow C/G$, and
$C_G(\sigma_1),\cdots,C_G(\sigma_r)$ the conjugacy classes whose elements generate the
stabilizers over the  branch points.

For  $\sigma\in G$ and $\chi\in \op{Irr}(G)$, we denote by
$l_{\sigma}(\chi)$ the multiplicity of the trivial character in the
restriction of $\chi$ to $\left<\sigma\right>$. Clearly,
$l_{\sigma}(\chi)\leq \chi(1)$, and the equality holds if and only
if $\sigma\in\op{Ker}(\chi)$.

 By \cite[Proposition~2]{Br87},
   for any nontrivial irreducible character $\chi$ of $G$,
\begin{equation*}\label{bro}
 h^1(C,\mathbb{C})^\chi=\chi(1)(2g(C/G)-2+r)-\sum_{j=1}^{r}l_{\sigma_j}(\chi),
\end{equation*} where
$h^1(C,\mathbb{C})^\chi =\dim \text{H}^i(C,\mathbb{C})^\chi$.

 In
particular,  if $g(C/G)=1$, then we have that,
\begin{equation}\label{cor-g=1}
\textrm{  $ h^1(C,\mathbb{C})^\chi\not=0$ if and only if
$\chi(\sigma_j)\not=\chi(1)$  for some $j$.}
\end{equation}

\subsection{The second cohomology of surfaces isogenous to a product of unmixed type}
Let $S=(C\times D)/G$ be a surface isogenous to a product of unmixed
type. Then the second cohomology of $S$ is
\begin{equation}\label{eq: cohomology}
 \begin{split}
\text{H}^2(S,\mathbb{C})& =\text{H}^2(C\times
D,\mathbb{C})^{\Delta_G}\\
 &=W \bigoplus  \left(\bigoplus_{\chi_1,\chi_2\in Irr(\Delta_G)}
  \text{H}^1(C,\mathbb{C})^{\chi_1}\otimes_\mathbb{C}
  \text{H}^1(D,\mathbb{C})^{\chi_2}\right)^{\Delta_G},
\end{split}\end{equation} where
$W=\text{H}^2(C,\mathbb{C})\otimes_\mathbb{C}\text{H}^0(D,\mathbb{C})\bigoplus
 \text{H}^0(C,\mathbb{C})\otimes_\mathbb{C}\text{H}^2(D,\mathbb{C})$ and $\Delta_G$ is the diagonal subgroup of $G\times G$. As a representation of $\Delta_G$, the irreducible consituents of $$\bigoplus_{\chi_1,\chi_2\in Irr(G)}
  \text{H}^1(C,\mathbb{C})^{\chi_1}\otimes_\mathbb{C}
  \text{H}^1(D,\mathbb{C})^{\chi_2}$$ all has the same character $\chi_1\chi_2$.
Hence the multiplicity of the trivial representation $\boldsymbol 1_{\Delta_G}$ in such a irreducible consituent is
\[
 \langle\chi_1\chi_2,\boldsymbol 1_{\Delta_G}\rangle_c = \langle\chi_1,\bar\chi_2\rangle_c =
\begin{cases}
 1, & \text{ if } \chi_2 = \bar\chi_1,\\
 0, & \text{ otherwise},
\end{cases}
\]
where $\langle\cdot\rangle_c$ is the inner product on the vector space of class functions on $G$. Therefore \\
$\left(\bigoplus_{\chi_1,\chi_2\in Irr(G)}
  \text{H}^1(C,\mathbb{C})^{\chi_1}\otimes_\mathbb{C}
  \text{H}^1(D,\mathbb{C})^{\chi_2}\right)^{\Delta_G} \neq 0$ if and only if $\chi_2=\bar\chi_1$,
and \eqref{eq: cohomology} becomes
\begin{equation}\label{eq: cohomology2}
 \text{H}^2(S,\mathbb{C}) = W \bigoplus  \left(\bigoplus_{\chi\in Irr(G)}
  \text{H}^1(C,\mathbb{C})^{\chi}\otimes_\mathbb{C}
  \text{H}^1(D,\mathbb{C})^{\bar\chi}\right)^{\Delta_G}.
\end{equation}

\begin{ex}\label{ex1}\rm {\bf Surfaces $S$ isogenous to a product of unmixed
type with   $\boldsymbol{\au\simeq\mathbb{Z}_2}$.}  Let $m, n, k, l$ be positive
integers. Let $G$ be one of the following groups:
 $$\mathbb{Z}_{2m}\oplus\mathbb{Z}_{2mn},\ \
 \mathbb{Z}_{2}\oplus\mathbb{Z}_{2m}\oplus\mathbb{Z}_{2mn}.$$
Let $\bar C$, $\bar D$ be elliptic curves. Let $\mathcal{V}:=(\alpha_1, \beta_1,
\overset{2k}{\overbrace{\gamma, \cdots, \gamma}})$,
$\mathcal{V}':=(\alpha_1', \beta_1',
\overset{2l}{\overbrace{\gamma', \cdots, \gamma'}}) $ be  generating
vectors of $G$ of type $(1; 2,\dots,2)$ with $\gamma\not=\gamma'$,
and $$\rho: C\to \bar C, \ \ \rho': D\to \bar D$$  the $G$-coverings
of smooth curves corresponding to $\mathcal{V}, \mathcal{V}'$,
respectively (cf.\,\ref{riemann}).

 For example, if $G=\mathbb{Z}_{2m}\oplus\mathbb{Z}_{2mn}$,
  we may take $\mathcal{V}=(\alpha, \beta,
\alpha^m, \cdots, \alpha^m)$, $\mathcal{V}'=(\alpha, \beta,
\beta^{mn}, \cdots, \beta^{mn})$, where  $\alpha:=(1, 0), \beta:=(0,
1)\in G$; if
$G=\mathbb{Z}_2\oplus\mathbb{Z}_{2m}\oplus\mathbb{Z}_{2mn}$,
  we may take $\mathcal{V}=(\mu, \nu,
\lambda, \cdots, \lambda)$, $\mathcal{V}'=(\mu, \nu, \lambda\mu^{m},
\cdots, \lambda\mu^{m})$, where  $\lambda:=(1, 0, 0), \mu:=(0, 1, 0)
,\nu:=(0,0,1)\in G$.

Let $G$  act diagonally on $C\times D$. Note that  the stabilizer of
each point lying  over any branch point of $\rho$ is $\left<
\gamma\right>$, and that of $\rho'$ is $\left< \gamma'\right>$ (cf.\,\ref{riemann}).  Since $\left< \gamma\right>\cap\left<
\gamma'\right>=1$ by assumption, we have that $G$ acts freely on
$C\times D$, and hence $S:=(C\times D)/G$ is a surface isogenous to
a product of curves.

By Hurwitz formula, we have $g(C)=2\delta m^2nk+1$ and $g(D)=2\delta
m^2nl+1$, where $\delta=1$ or $4$ depending  on
$G=\mathbb{Z}_{2m}\oplus\mathbb{Z}_{2mn}$ or not.
 So the numerical invariants of $S$ are as below:
\begin{align*}\textrm{ $p_g(S)=\delta m^2nkl+1$, $q(S)=2$ and
$K_S^2=8\delta m^2nkl$.}
\end{align*}

Let $I=\{\chi\in \op{Irr}(G)\ |\ \chi(\gamma)\not=1 \ \textrm{and}
 \ \bar\chi(\gamma')\not=1\}$.
 By \eqref{cor-g=1} and \eqref{eq: cohomology2}, we have
\begin{equation}\label{eq1}
 \begin{split}
\text{H}^2(S,\mathbb{C})=W   \bigoplus_{\chi\in I}
 \text{H}^1(C,\mathbb{C})^\chi\otimes_\mathbb{C}
  \text{H}^1(D,\mathbb{C})^{\bar\chi},
\end{split}
\end{equation}
with
$W=\text{H}^2(C,\mathbb{C})\otimes_\mathbb{C}\text{H}^0(D,\mathbb{C})\bigoplus
 \text{H}^0(C,\mathbb{C})\otimes_\mathbb{C}\text{H}^2(D,\mathbb{C})$.

Since $\gamma$ (resp. $\gamma'$) is of order two, it induces
$-\op{id}$ on $\text{H}^1(C,\mathbb{C})^\chi$ (resp.
$\text{H}^1(D,\mathbb{C})^{\bar\chi}$) for all $\chi\in I$. So
$(\gamma , \gamma')$ induces identity on the right-side hand of
(\ref{eq1}).

Let $\sigma$ be the automorphism  of $S$ induced by $(\gamma ,
\gamma')\in  \op{Aut}(C)\times \op{Aut}(D)\subseteq\op{Aut}(C\times
D)$. Then $\sigma$ is an involution of $S$ and it acts trivially on
$H^2(S, \mb Q)$ and hence on $H^*(S, \mb Q)$.
\end{ex}

\begin{rmk}\rm Surfaces in Example~\ref{ex1}
are  rigidified by Proposition~\ref{bn} below. It is not known
whether they are cohomologically rigidified.
\end{rmk}
\begin{prop}\label{bn}Let $S$ be a smooth projective surface.
Assume that the universal cover  of $S$ is a bounded domain in
$\CC^2$. Then $S$ is rigidified.
\end{prop}
\begin{pf} Otherwise, there is an automorphism $\sigma\in \op{Aut}(X)\cap
\op{Diff}^0(X)$, of prime order. The assumption implies $S$ is of
general type; in particular $\chi(\omega_S)>0$. So we have
$S^\sigma\not=\emptyset$ by the proof of Theorem~\ref{mad}. Thus
$\sigma$ and $\op{id}_S$ are  homotopic automorphisms which agree at
 $S^\sigma$. Since the universal cover of $S$ is a product of two unit disks in $\CC$, a bounded domain, it follows that $\sigma=\op{id}_S$ by \cite[Theorem~3.6]{BN67} --- a contradiction.
\qed\end{pf}
\begin{thm}\label{non-abelian}
 Let $S=(C\times D)/G$ be a surface isogenous to a product of unmixed type
 with $q(S)\geq2$.
 If $S$ is  rationally cohomologically rigidified, then $S$ is as in Example~\ref{ex1}.
\end{thm}

Before proving Theorem~\ref{non-abelian}, we show the following
preparatory results.
\begin{prop}\label{proposition: automorphism1}
Let $S=(C\times D)/G$ be a surface isogenous to a product (with
minimal realization).  Denote by  $\Delta_G$ the diagonal of
$G\times G$.  Then
\[
 \op{Aut}(S)=N(\Delta_G)/\Delta_G,
\]
where $N(\Delta_G)$ is the normalizer of $\Delta_G$ in
$\op{Aut}(C\times D)$.
\end{prop}
\begin{pf}
 For each
$\sigma\in \op{Aut}(S)$, there is an automorphism $\tilde{\sigma}\in
Aut(C\times D)$ such that
\[
\begin{CD}
 C\times D @>\tilde{\sigma}>> C\times D \\
 @V\pi VV                          @VV\pi V \\
 S         @>\sigma>>       S
\end{CD}
\]
is commutative, where $\pi$ is the quotient map. The existence of
such a lift $\tilde{\sigma}$ of $\sigma$ follows  simply from
 the uniqueness of minimal realization of $S$.

 On the other hand, given $\tilde{\sigma}\in
\op{Aut}(C\times D)$, $\tilde{\sigma}$ descends to an automorphism
$\sigma\in \op{Aut}(S)$ if and only if it is in the normalizer
$N(\Delta_G)$ of $\Delta_G$ in $\op{Aut}(C\times D)$. Hence we have
a surjective homomorphism of groups $N(\Delta_G)\rightarrow
\op{Aut}(S)$ and its kernel is easily seen to be $\Delta_G$. So
$\op{Aut}(S)=N(\Delta_G)/\Delta_G$.\qed
\end{pf}

\begin{prop}\label{proposition: aut0}
Let $S$ be as in Proposition~\ref{proposition: automorphism1}.  If
$S$ is of unmixed type, then
\[
\op{Aut}_0(S) \subseteq (G\times G)\cap N(\Delta_G)/\Delta_G.
\]
\end{prop}
\begin{pf}
For each $\sigma\in \op{Aut}_0(S)$, let $\tilde\sigma\in
\op{Aut}(C\times D)$ be its lift as in the proof of Proposition~\ref{proposition: automorphism1}. By the proof of Lemma~\ref{3.1},
$\sigma$ preserves the two induced  fibrations $\pi_1\colon S
\rightarrow C/G$ and $\pi_2\colon S \rightarrow D/G$, and it induces
identity  on their bases $C/G$ and $D/G$.  Hence   $\tilde \sigma$
does not interchange the factors of $C\times D$. By \cite[Rigidity
Lemma~3.8]{Cat00}), there are automorphisms  $\sigma_C$ and
$\sigma_D$, of $C$ and $D$, respectively,  such that $\tilde
\sigma=(\sigma_C,\sigma_D)$. Since  $\tilde \sigma$ induces identity
on  bases of $\pi_1$ and $\pi_2$, we have $\sigma_C,\sigma_D\in G$.

On the other hand, By Proposition~\ref{proposition: automorphism1},
we have $\tilde \sigma \in N(\Delta_G)$. So $\tilde\sigma\in
(G\times G)\cap N(\Delta_G)$, and
 $\sigma=\tilde\sigma\mod\Delta_G\in (G\times G)\cap N(\Delta_G)/\Delta_G$.\qed
\end{pf}
\begin{rmk}\rm\label{remark: centralizer}
Let $Z_G$ be the center of $G$. Then $(G\times G)\cap N(\Delta_G)$
is generated by $Z_G\times\{1\}$  and $\Delta_G$, and the map
$$Z_G\to (G\times G)\cap N(\Delta_G)/\Delta_G,\ \sigma\mapsto (\sigma, 1)
 \mod\Delta_G$$
 is an isomorphism of groups.
In what follows, we regard $Z_G$ as a subgroup of $\op{Aut}(S)$ under
 such an isomorphism.

 So by Proposition~\ref{proposition: aut0}, $\au$ is (isomorphic to) a subgroup of
$Z_G$; in particular, if $G$ is centerless, then $\au$ is trivial.
\end{rmk}

\begin{lem}\label{eq2} Let $S$ be as in Proposition \ref{proposition: automorphism1}.
 If
$S$ is of unmixed type, then for each
 $\sigma\in Z_G\subseteq \op{Aut}(S)$, we have that, $\sigma\notin \au$
   if and only if there is an
irreducible $\chi\in \op{Irr}(G)$ such that $\sigma\notin
\op{Ker}(\chi)$,
 $\op{H}^1(C,\mathbb{C})^\chi\neq 0$,    and
  $\op{H}^1(D,\mathbb{C})^{\bar\chi}\neq 0$.
\end{lem}
\begin{pf}
  Since  for each
$\sigma\in Z_G\subseteq \op{Aut}(S)$, $(\sigma,1)\in
\op{Aut}(C\times D)$ is a lift of $\sigma$ (cf. Remark \ref{remark:
centralizer}), the lemma follows from the fact that, for each
$\sigma\in \op{Aut}(S)$,  the quotient map $C\times D\to S$ induces
an isomorphism between the action of $\sigma$ on
$\text{H}^2(S,\mathbb{C})$ and that of its lift $\tilde{\sigma}$ on
the right-hand side of (\ref{eq: cohomology2}). \qed\end{pf}

\begin{lem}\label{induced}
 Let  $H$ be a subgroup of a finite group $G$,
 and $\chi$  an irreducible character of $H$.
   Let $H'\subseteq H$ be  a subset such that $H'\cap\op{Ker}(\chi)=\emptyset$.
   Then
\begin{enumerate}\item[(i)]
 for any irreducible constituent
    $\varphi$ of $\chi^G$, $H'\cap\op{Ker}(\varphi)=\emptyset$;
\item [(ii)]
  if moreover $\chi^G(g)=0$ for some $g\in G$,
then there is an irreducible constituent
    $\varphi'$ of $\chi^G$ such that
     $(\{g\}\cup H')\cap\op{Ker}(\varphi')=\emptyset$.
    \end{enumerate}
\end{lem}
\begin{pf} For any irreducible constituent
    $\varphi$ of $\chi^G$,
Frobenius reciprocity theorem gives
 $(\varphi,\chi^G) = (\varphi|_H,\chi)$.
 Hence the multiplicity of $\chi$ in $\varphi|_H$ is the same as
  that of
  $\varphi$ in $\chi^G$. In particular
  $\chi$ is a constituent of $\varphi|_H$ and
  $\op{Ker}(\varphi)\cap H \subseteq \op{Ker}(\chi)$. So (i)
  follows.

If  $\chi^G(g)=0$ for some $g\in G$, then there exists an
irreducible constituent
    $\varphi'$ such that   $g\notin \op{Ker}(\varphi')$. Combining
    this with (i), we obtain (ii).
\qed\end{pf}

\subsection{Proof of Theorem~\ref{non-abelian}}
  Suppose that $\op{Aut}_0(S)$ is not trivial.   By
Proposition~\ref{mtype}, we may assume that $S$ is of
 unmixed type.
Consider the induced fibrations $S\to C/G$ and $S\to
 D/G$. By Lemma~\ref{3.1}, note that $q(S)=g(C/G)+g(D/G)$,
  we have that $g(C/G)=g(D/G)=1$.

 Let $\mathcal{U}:=(a,b,\sigma_1,\dots,\sigma_r)$
  (resp. $\mathcal{U}':=(c,d,\tau_1,\dots,\tau_s)$) be the generating vector of $G$
  for the branch covering $C\rightarrow C/G$ (resp. $D\rightarrow
  D/G$) (cf.\,\ref{rie}). Denote by $m_i$, $n_j$ the order of $\sigma_i$, $\tau_j$,
  respectively.

Let $ \Sigma_1 = \cup_{g\in G}\cup_{1\leq i\leq r}
   \left<g\sigma_i g^{-1}\right>$
   (resp. $\Sigma_2 = \cup_{g\in G}\cup_{1\leq j\leq s}
   \left<g\tau_j g^{-1}\right>$).
Since the action of $G$ on $C\times D$ is free, we have
$\Sigma_1\cap \Sigma_2 =\{1\}$.

By Theorem~\ref{main}, $\au$ is of order two. Let  $\sigma\in Z_G$
such that $\sigma$ is the generator of $\au$ (cf.\,Remark~\ref{remark: centralizer}).
 By Lemma \ref{eq2} and  (\ref{cor-g=1}), we
have that
\begin{equation}\label{not=1}
 \textrm{ for any $\chi\in I$, $\sigma\in \op{Ker}(\chi)$,}
\end{equation}
where $I$ is the set of irreducible characters $\chi$ of $G$  such
that
   $\sigma_i, \tau_j\notin \op{Ker}(\chi)$
 for some $i,j$.

\subsubsection{Claim: $G$ is  abelian.}
\subsubsection{Proof of Claim.}
If $G$ is not abelian, we will get a contradiction by finding  an
irreducible character $\chi\in I$ such that $\sigma\notin
\op{Ker}(\chi)$. By Lemma~\ref{induced}, it is enough to find a
subgroup $H$ of $G$ and an irreducible character $\chi$ of $H$, such
that $\sigma,\sigma_i,\tau_j\in H$ and $\sigma,\sigma_i,\tau_j\notin
\op{Ker}(\chi)$ for some $i, j$.

 For each $1\leq i\leq r$ and $1\leq j\leq s$, let
$G_{ij}$ be the subgroup of $G$ generated by $\sigma_i$ and
$\tau_j$, $\varphi_i$ the linear character of the cyclic group
$\left<\sigma_i\right>$ such that $\varphi_i(\sigma_i)=\xi$, where
$\xi$ is a primitive $m_i$-th root, and  $\varphi_i^{G_{ij}}$   the
induced character from $\varphi_i$. Since
$\Sigma_1\cap\Sigma_2=\emptyset$, we have
\[
 \varphi_i^{G_{ij}}(\tau_j) =0
\]
for all $1\leq i\leq r$ and $1\leq j\leq s$. By Lemma~\ref{induced},
   there is an irreducible character $\chi_{ij}$ of
$\varphi_i^{G_{ij}}$ such that
\begin{equation}\label{non-ab}
 \sigma_i,\tau_j \notin \op{Ker}(\chi_{ij}).
\end{equation}

Similarly, starting with a primitive linear character of
$\left<\tau_j\right>$, we can construct  a  character $\chi_{ji}'$
of $G_{ij}$
 such that
\[
 \sigma_i,\tau_j \notin \op{Ker}(\chi_{ji}').
\]

\textbf{Step 1.} First we assume  $\sigma\notin G_{ij}$ for some $i, j$. Since
$\sigma$ is of order two and $\sigma\in Z_G$, we have
$\left<\sigma\right> \cap G_{ij} =\{1\}$, and  for any $g\in G$,
$g\sigma g^{-1}=\sigma \notin G_{ij}$. So by the definition of
induced character, we have
$\chi_{ij}^{\left<G_{ij},\sigma\right>}(\sigma)=0$. By Lemma~\ref{induced}, there is an irreducible character $\tilde\chi_{ij}$
of $\left<G_{ij},\sigma\right>$ such that
$\sigma,\sigma_i,\tau_j\notin \op{Ker}(\tilde\chi_{ij})$.

\textbf{Step 2.} Next, suppose $\sigma\in G_{ij}$ for any $1\leq i\leq r$ and $1\leq
j\leq s$.

If $\sigma\in \left<\sigma_i\right>$ for some $1\leq i\leq r$, then
$\varphi_i(\sigma)\not=1$. By Lemma~\ref{induced} and
(\ref{non-ab}), $\sigma,\sigma_i,\tau_j \notin \op{Ker}(\chi_{ij})$.

Similarly, If $\sigma\in \left<\tau_j\right>$ for some $1\leq j\leq
s$, then $\sigma,\sigma_i,\tau_j \notin \op{Ker}(\chi_{ji}')$.

\textbf{Step 3.} Now we  can  assume that, ($\ast$) $\sigma\in G_{ij}$ but
$\sigma\notin\left<\sigma_i\right>$ and
$\sigma\notin\left<\tau_j\right>$ for all $1\leq i\leq r$ and $1\leq
j\leq s$.

Since $\sigma$ is in the center of $G$, we have
$\left<\sigma,\sigma_i\right> \cong \ZZ_2\oplus\ZZ_{m_i}$. Let
$\psi_i$ be the linear character of $\left<\sigma,\sigma_i\right>$
such that
\[
 \psi_i(\sigma)=-1,  \ \  \psi_i(\sigma_i) =\xi,
\]
where $\xi$ is a primitive $m_i$-th root. Let $_j\psi_i$ be the
induced character of $\psi_i$ on the  group
$\left<\sigma,\sigma_i,\tau_j\right>$.

\textbf{Step 3.1.} If  $C_G(\tau_j)\cap
\left<\sigma,\sigma_i\right>=\emptyset$ for some $1\leq i\leq r$ and
$1\leq j\leq s$, where $C_G(\tau_j)$ is the conjugate class of
$\tau_j$ in $G$,  then by the definition of induced character,  $
_j\psi_i(\tau_j) =0$. By Lemma~\ref{induced}, there is a constituent
$\psi$ of $_j\psi_i$ such that $\sigma,\sigma_i, \tau_j\notin
\op{Ker}(\psi)$.

\textbf{Step 3.2.} Next,  we assume additionally $C_G(\tau_j)\cap
\left<\sigma,\sigma_i\right> \neq \emptyset$ for all $1\leq i\leq r$
and $1\leq j\leq s$.  Then for each  $1\leq i\leq r$ and $1\leq
j\leq s$, there is an element $g_{ij}\in G$ such that $g_{ij}\tau_j
g_{ij}^{-1}\in \left<\sigma,\sigma_i\right>$.

If $m_i\geq3$ for some  $1\leq i\leq r$, then it is easy to find
 a linear character $\chi$ of $\left<\sigma,\sigma_i\right>$
such that $\sigma,\sigma_i,g_{ij}\tau_j g_{ij}^{-1}\notin
\op{Ker}(\chi)$ and hence $\sigma,\sigma_i,\tau_j \notin
\op{Ker}(\chi)$. (Since  $H_i:= \left<\sigma,\sigma_i\right>$ is
isomorphic to $\ZZ_{2} \oplus \ZZ_{m_i}$,
  we can find characters  $\phi, \phi'\in H_i$ such
that
\[
 \phi(\sigma)=1, \ \phi(\sigma_i)=\xi_i,\
 \phi'(\sigma_i)=-1,  \ \phi'(\sigma_i)=1,
\]
where $\xi_i$  is  a root of unity of order $m_i$. Write
$\tau_j':=g_{ij}\tau_j g_{ij}^{-1}=\sigma^a\sigma_i^b$ for some $a,
b\geq0$. Since $\tau_j'\not=\sigma$, we have $b\not=0$. Let
$\phi_k=\phi^k\phi'$ for $k=1, 2$. Then
$\phi_k(\tau_j')=(-1)^a\xi_i^{bk}$. So among $\phi_1$ and $\phi_2$,
there is at least one character, say $\phi_1$, such that
$\phi_1(\tau_j')\not=1$.)

Similarly, if $n_j\geq3$ for some  $1\leq j\leq s$, then we can find
 a linear character $\chi$ of $\left<\sigma,\tau_j\right>$
such that  $\sigma,\sigma_i,\tau_j \notin \op{Ker}(\chi)$.

\textbf{Step 3.3.} Finally we may assume further  $m_i = n_j =2$  for all $1\leq i\leq
r$ and $1\leq j \leq s$. This implies that $G_{ij}=D_{2k_{ij}}$ for
some $k_{ij}\geq2$ since any finite group generated by two elements
of order two is dihedral.

If $k_{ij}>2$ for some $i$ and $j$, then it is well known that there
is a faithful ($2$-dimensional) representation $\rho$ of $G_{ij}$.
Let $\chi$ be the corresponding character of $\rho$.  Since $\rho$
is faithful,   $\sigma,\sigma_i, \tau_j\notin \op{Ker}(\chi)$.

If $k_{ij}=2$ for all $i$ and $j$, the assumption ($\ast$) above
implies that $\sigma=\sigma_i\tau_j$  for all $1\leq i\leq r$ and
$1\leq j \leq s$. So $\sigma_1=\dots=\sigma_r$,  $\tau_1=\dots
=\tau_s$ and $\sigma_1\not=\tau_1$.

\textbf{Step 4.} We show that  $\sigma_1$, $\tau_1$ are in the center of $G$.

Note that  $G_{ij}$ in the proof above can be replaced by
$G_{ij}'=\left<\sigma_i', \tau_j' \right>$ for any $ \sigma_i'\in
C_G(\sigma_i), \tau_j'\in C_G(\tau_j)$, since the characters of $G$
do not distinguish conjugate elements. So  we can assume much more,
namely, $\left<\sigma_i',\tau_j'\right>\cong \mathbb Z^{\oplus 2}$,
$\sigma\in \left<\sigma_i',\tau_j'\right>$,  $\sigma\notin
\left<\sigma_i'\right> $ or $\left<\tau_j'\right>$ for any $i,j$,
$\sigma_i'\in C_G(\sigma_i)$ and $\tau_j'\in C_G(\tau_j)$. Under
these assumptions, we have
\[\sigma\in
A:=\left<\sigma_1',\tau_1\right> \cong \mathbb Z_2^{\oplus 2}
\]
for any $\sigma_1'\in C_G(\sigma_1)$. This implies that $A$ is
generated by $\sigma$ and $\tau_1$, and hence it is generated by
$\sigma_1$ and $\tau_1$ since $\sigma=\sigma_1\tau_1$. So we  have
 $\sigma_1'=\sigma_1$ and hence  $C_G(\sigma_1)=\left<\sigma_1\right>$.

   Since $\sigma\in Z_G$ and $\sigma=\sigma_1\tau_1$,
    we have that $\tau_1$ is in the center of $G$.

Now by the definition of a generator vector, $G$ is generated by  $a$, $b$ and $\tau_1$, and
$aba^{-1}b^{-1}\sigma_1^r=1$.

 This implies that the
 commutator subgroup $G'$ is contained in $\left<\sigma_1\right>$.
Similarly $G'$ is a subgroup of $\left<\tau_1\right>$. Since
$\left<\sigma_i\right>\cap \left<\tau_j\right> = \{1\}$, $G'$ is
trivial and $G$ is abelian --- a contradiction. This finishes the proof
of the claim.\qed
\\

Now we may assume that $G$ is abelian. By the proof
of the claim, we have that $\sigma_1=\dots=\sigma_r$,  $\tau_1=\dots
=\tau_s$ and $\sigma=\sigma_1\tau_1$. Thus $G$ can be generated by
three elements, namely $a,b$ and $\sigma_1$ (or $c,d$ and $\tau_1$).
By the structure theorem of finitely generated abelian groups, we
may write $G=\ZZ_{d_1}\oplus\ZZ_{d_2}\oplus\ZZ_{d_3}$ with $d_1\ |\
d_2\ |\ d_3$. Since $G$ has at least two elements of order 2, both
$d_2$ and $d_3$ are even. If $d_1=1$, then  $G\cong\ZZ_{2m} \oplus
\ZZ_{2mn}$ for some positive integers $m, n$.

If $d_1\geq 2$, then $G$ needs three generators, one of which is
$\sigma_1$ or $\tau_1$. Since $\sigma_1$ and $\tau_1$ have order
$2$, we see that $d_1=2$. Hence in this case $G\cong \ZZ_{2}\oplus
\ZZ_{2m}\oplus\ZZ_{2mn}$ for some positive integers $m,  n$.

Since $aba^{-1}b^{-1}\sigma_1^r=1$ and  $cdc^{-1}d^{-1}\tau_1^s=1$
in $G$, we have that both $r$ and $s$ are even. So $S$ is as in
Example~\ref{ex1} with $\mathcal{V}=\mathcal{U}$ and
$\mathcal{V}'=\mathcal{U}'$.

This completes the proof of Theorem~\ref{non-abelian}. \qed

\end{document}